\documentclass{amsart}
\usepackage{amssymb}

\usepackage{graphicx}
\usepackage{amscd}

\newtheorem{theorem}{Theorem}
\theoremstyle{plain}

\newtheorem{definition}{Definition}

\newtheorem{lemma}{Lemma}

\numberwithin{equation}{section}
\numberwithin{theorem}{section}
\numberwithin{lemma}{section}
\numberwithin{proposition}{section}
\numberwithin{corollary}{section}
\numberwithin{remark}{section}
\numberwithin{example}{section}
\numberwithin{definition}{section}

\input{tcilatex}

\begin{document}
\title[Poles of Local Zeta Functions]{On the poles of Igusa's Local Zeta Function for algebraic sets}
\author{W. A. Zuniga-Galindo}
\address{Barry University, Department of Mathematics and Computer Science\\
11300 N.E. Second Avenue, Miami Shores, Florida 33161, USA}
\email{wzuniga@mail.barry.edu}

\begin{abstract}
Let $K$ be a $p-$adic field, and $Z_{\Phi }(s,f)$, $\,\,\,s\in \mathbb{C}$,
with $Re(s)>0$, the Igusa local zeta function associated to $%
f(x)=(f_{1}(x),..,f_{l}(x))$ $\in \left[ K\left( x_{1},..,x_{n}\right) 
\right] ^{l}$, and $\Phi $ a Schwartz-Bruhat function. The aim of this paper
is to describe explicitly the poles of the meromorphic continuation of $%
Z_{\Phi }(s,f)$. Using resolution of singularities  is possible to express $%
Z_{\Phi }(s,f)$\ as a finite sum \ of $p-$adic monomial integrals. These
monomial \ integrals are computed explicitly by using techniques of toroidal
geometry. In this way, an explicit list for the candidates to poles of $%
Z_{\Phi }(s,f)$ \ is obtained. \ 
\end{abstract}

\maketitle

\section{\protect\bigskip Introduction}

Let $K$ be a $p-$adic field, i.e. $[K:\mathbb{Q}_{p}]<\infty $. Let $%
\mathcal{O}_{K}$\ be the valuation ring of $K$, $\mathcal{P}_{K}$ the
maximal ideal of $\mathcal{O}_{K}$, and $\overline{K}=\mathcal{O}_{K}/%
\mathcal{P}_{K}$ \ the residue field of $K$. \ The cardinality of the
residue field of $K$ is denoted by $q$, thus $\overline{K}=\mathbb{F}_{q}$.
For $z\in K$, $v(z)\in \mathbb{Z}\cup \{+\infty \}$ \ denotes the valuation
of $z$, $\left| z\right| _{K}=q^{-v(z)}$, and $ac$ $z=z\pi ^{-v(z)}$ where $%
\pi $\ is a fixed uniformizing parameter for $\mathcal{O}_{K}$. For $%
x=(x_{1},...,x_{l})\in K^{l}$, we define $\left\| x\right\|
_{K}:=\max_{1\leq i\leq l}\left| x_{i}\right| _{K}$.

Let $f_{i}(x)$ $\in K[x]$, $x=(x_{1},...,x_{n})$, be a non-constant
polynomial for $i=1,...,l$. Let $\Phi :K^{n}\rightarrow \mathbb{C}$ be a
Schwartz-Bruhat function, i.e. \ a locally constant function with compact
support. The Igusa local zeta function attached to $f(x):=$ $\left(
f_{1}(x),...,f_{l}(x)\right) $ is defined as

\begin{equation}
Z_{\Phi }(s,f)=\int\limits_{K^{n}}\Phi (x)\left\| f(x)\right\| _{K}^{s}\mid
dx\mid ,\,\,\,\,\,\mathbb{\,}\noindent \,\,\,s\in \mathbb{C},  \label{1}
\end{equation}
for $Re(s)>0$, where $\mid dx\mid $ denotes the Haar measure on $K^{n}$ so
normalized that $\mathcal{O}_{K}^{n}$ has measure $1$. Since the numerical
data of a resolution \ of singularities of $\cup _{i=1}^{l}f_{i}^{-1}(0)$
are directly related to the meromorphic continuation of $Z_{\Phi }(s,f)$, \
we say that $Z_{\Phi }(s,f)$ is the local zeta function \ associated to $%
\cup _{i=1}^{l}f_{i}^{-1}(0)$. Alternatively, we may say \ that $Z_{\Phi
}(s,f)$ is the local zeta function \ associated to the polynomial mapping $%
f:K^{l}\rightarrow K^{n}$.

The functions $Z_{\Phi }(s,f)$ were introduced by Weil \cite{W} and their
basic properties for general $f$ and $l=1$ were first studied by Igusa.
Using resolution of singularities Igusa proved that $Z_{\Phi }(s,f)$ admits
a meromorphic continuation to the complex plane as a rational function of $%
q^{-s}$ \cite{I1}, \cite[Theorem 8.2.1]{I2}. Igusa's proof was generalized \
by Meuser to the case $l\geqslant 1$ \cite[Theorem 1]{M2}. Using $p-$adic
cell decomposition, Denef gave a completely different proof of the
rationality of $Z_{\Phi }(s,f)$, with $l\geqslant 1$ \cite{D1}.

\bigskip We write $Z(s,f)$ \ \ when $\Phi $ is the characteristic function
of $\mathcal{O}_{K}^{n}$.\ Suppose that $f_{i}(x)$, $i=1,...,l$, have
coefficients in $\mathcal{O}_{K}$. Let $N_{j}(f)$ be the number of solutions
of $f_{i}(x)\equiv 0\,\ $mod$\,\ \mathcal{P}_{K}^{j},\ i=1,...,l$, in $%
\left( \mathcal{O}_{K}/\mathcal{P}_{K}^{j}\right) ^{n}$,$\,\ $and put \ $%
\,P(t,f):=\sum_{j=0}^{\infty }N_{j}(f)(q^{-n}t)^{j}$. The Poincar\'{e}
series $P(t,f)$ is related to $Z(s,f)$ by the formula \ 

\begin{equation}
P(t,f)=\frac{1-tZ(s,f)}{1-t},\text{ \ \ }t=q^{-s},  \label{2}
\end{equation}
\cite[Theorem 2]{M2}. Thus the rationality of \ $Z(s,f)$ implies the
rationality of \ $P(t,f)$. In the case $l=1$ the rationality of $P(t,f)$ was
conjectured by S. E. Borewicz and I. R. \v{S}afarevi\v{c} \cite[page 63]{BS}%
. The Borewicz-\v{S}afarevi\v{c} conjecture \ is a consequence of the
rationality of $Z(s,f)$ proved by Igusa, Meuser, and Denef. \ 

\bigskip

A basic problem is to determine the poles of the meromorphic continuation of 
$Z_{\Phi }(s,f)$ in $\func{Re}(s)<0$. In the case $l=1$, the basic strategy
is to take a resolution of singularities $h:X_{K}\rightarrow A_{K}^{n}$ of $%
\ f$ $^{-1}(0)$ and study the numerical data $\{(N_{i},v_{i})\}$ of the
resolution, where $N_{i}$ is the multiplicity of $f$ along the exceptional
divisor $D_{i}$, and $\nu _{i}-1$ is the multiplicity of $h^{\ast
}(\bigwedge\limits_{1\leq j\leq n}dx_{j})$ along the divisor $D_{i}$. Igusa
showed that the set of ratios $\left\{ -\frac{v_{i}}{N_{i}}\right\} \cup
\left\{ -1\right\} $ contains the real parts of the poles of $Z_{\Phi }(s,f)$
\cite{I1}, \cite[Theorem 8.2.1]{I2}. Many examples show that almost all
these fractions do not correspond to poles of $Z_{\Phi }(s,f)$. In the case
of polynomials in two variables, as a consequence of the works of Igusa,
Strauss, Meuser and Veys, there is a complete solution of this problem \ 
\cite{I1a}, \cite{S}, \cite{M1}, \cite{V1}, \cite{V2}. For general
polynomials the problem of determination of the poles of $\ Z_{\Phi }(s,f)$
is still open. There exists a generic class of polynomials called
non-degenerate with respect to its Newton polyhedron for which is possible
to give a small set of candidates for the poles of $Z_{\Phi }(s,f)$. The
poles of the local zeta functions attached to non-degenerate polynomials can
be described in terms of Newton polyhedra. The case of two variables was
studied by Lichtin and Meuser \cite{LM}. In \cite{D2}, Denef gave a
procedure based on monomial changes of variables to determine a small set of
candidates for the poles of $Z_{\Phi }(s,f)$ in terms of the Newton
polyhedron of $f$. This result was obtained by the author, using an approach
based on the $p-$adic stationary phase formula and N\'{e}ron $p-$%
desingularization, for polynomials with coefficients in a non-archimedean
local field of arbitrary characteristic \cite{Z2}, (see also \cite{DH}, \cite
{Z1}, \cite{Z3}).

\bigskip The local zeta functions $Z_{\Phi }(s,f)$, \ with $l\geqslant 1$,
are a generalization of the Igusa local zeta functions \ associated to
hypersurfaces. Another generalizations are due to Denef and Loeser. \ In 
\cite{L} Loeser \ introduced a \ local zeta function in $l$ complex
variables associated to an algebraic set of the \ form $D_{K}$. Recently,
Denef and Loeser have introduced \ \ a completely new class of zeta
functions called motivic Igusa zeta functions \cite{DL1} that include as
particular cases the Igusa local zeta functions, with $l=1$, and the
topological zeta functions \cite{DL2}.

\ The first result of this paper (see theorem \ref{mainth}) provides an
explicit description of the poles of \ $Z_{\Phi }(s,f)$, \ with $l\geqslant
1 $, in terms of a resolution of singularities of $f$. Using resolution of
singularities \ it is possible to express $Z_{\Phi }(s,f)$\ as a finite sum
\ of $p-$adic monomial integrals. We compute explicitly these monomial
integrals by using toroidal geometry (see section 3). In this way we are
able to give an explicit list for the candidates to poles of $Z_{\Phi }(s,f)$
\ in terms of \ a list of Newton polyhedra constructed from the numerical
data associated to a resolution of \ singularities of the divisor $\cup
_{i=1}^{l}f_{i}^{-1}(0)$. The second result \ provides an explicit
description of the poles of a twisted local zeta function $Z_{\Phi }(s,\chi
,f)$ in terms of a resolution of singularities of $f$ (see theorem \ref
{theorem2}). Finally, we give \ a family of algebraic sets, that we called
monomial algebraic sets, for which the corresponding local zeta functions
can be expressed in terms of Newton polyhedra (see theorem \ref{th}).

\section{Preliminaries}

\subsection{Numerical data}

We put $D_{K}:=\cup _{i=1}^{l}f_{i}^{-1}(0)$. By Hironaka's resolution
theorem \cite{H} applied to $D_{K}$, there exist an $n-$dimensional $K-$%
analytic manifold $X_{K}$ and proper $K-$analytic map \ $h:X_{K}\rightarrow
A_{K}^{n}$ with the following properties:

(1) $h^{-1}\left( D_{K}\right) $ is a divisor with normal crossings;

(2) the restriction of $h$ to \ $X_{K}$ $\setminus h^{-1}\left( D_{K}\right) 
$ is an isomorphism onto its image.

We denote by $E_{i}$, $i\in T$, the irreducible components of the divisor $%
h^{-1}\left( D_{K}\right) $. At every point $b$ of $X_{K}$ \ if $E_{1}$, $%
E_{2}$,$...$, $E_{r}$, $r\leq n$, are all the components containing $b$ with
respective local equations $y_{1}$,$...$,$y_{r}$ around $b$, \ then there
exist local coordinates of $X_{K}$ around $b$ of \ the form $(\
y_{1},...,y_{r},y_{r+1},...,y_{n})$ \ such that 
\begin{equation}
f_{i}\circ h=\epsilon _{i}\prod\limits_{1\leq j\leq r}y_{j}^{N_{i,j}}\text{, 
}i=1,...,l,  \label{3}
\end{equation}
\begin{equation}
h^{\ast }(\bigwedge\limits_{1\leq j\leq n}dx_{j})=\left( \eta
\prod\limits_{1\leq j\leq r}y_{j}^{v_{j}-1}\right) \bigwedge\limits_{1\leq
j\leq n}dy_{j},  \label{4}
\end{equation}
on some neighborhood of $b$, in which $\epsilon _{i}$, $i=1,...,l$, $\eta $\
are units of the local ring $\mathcal{O}_{b}$ of \ $X_{K}$ at $b$.

We set 
\begin{eqnarray}
N_{i} &:&=(N_{i,1},...,N_{i,r},0,..,0)\in \mathbb{N}^{n}, \\
N_{b} &:&=\left\{ N_{1},...,N_{l}\right\} ,  \notag
\end{eqnarray}
and 
\begin{equation}
v_{b}:=(v_{1}-1,...,v_{r}-1,0,...,0)\in \mathbb{N}^{n}.
\end{equation}
We define the \textit{numerical data} of $\ f=(f_{1},...,f_{l})$ as the set
of all pairs $(N_{b},v_{b})$ satisfying $N_{i}\neq 0$, $i=1,...,l$, and $%
b\in $ $X_{K}$.

We note that if \ $(N_{b}^{\prime },v_{b}^{\prime })$ is a numerical datum
at $b\in X_{K\text{ }}$ obtained by using a second local coordinate system,
then it may be verified that \ $N_{b}=\ N_{b}^{\prime }$, $%
v_{b}=v_{b}^{\prime }$, after a permutation of indices, i.e. after a
renaming of \ the local coordinates \ around $b$. As we shall see later on
our description of the poles of \ $Z_{\Phi }(s,f)$ does not depend on the
local \ coordinate systems used in the computation of the numerical data of $%
\ f$.

\subsection{Newton polyhedra}

Let $\mathbb{R}_{+}:=\{x\in \mathbb{R}\mid x\geqslant 0\}$, and $S$ a
non-empty subset of $\mathbb{N}^{n}$. The \textit{Newton polyhedron} $\Gamma
(S)$ \ associated to $S$ is defined as the convex hull \ in $\mathbb{R}%
_{+}^{n}$ of the set $\bigcup\limits_{m\in S}\left( m+\mathbb{R}%
_{+}^{n}\right) $. \ A \textit{facet} is a face of \ $\Gamma (S)$ of
dimension $n-1$.

We denote by $\left\langle \text{ }.\text{ },.\text{ }\right\rangle $ the
usual inner product of $\mathbb{R}^{n}$. For $a\in \mathbb{R}_{+}^{n}$, we
define 
\begin{equation}
m(a):=\inf_{x\in \Gamma (S)}\left\langle a,x\right\rangle .
\end{equation}

Given $a\in \mathbb{R}_{+}^{n}$, the \ \textit{first meet locus} $F(a)$ of $%
a $ is defined as 
\begin{equation}
F(a):=\{x\in \Gamma (S)\mid \left\langle a,x\right\rangle =m(a)\}.
\end{equation}
The first meet locus \ is a face of $\Gamma (S)$. Moreover, if $a\neq 0$, $%
F(a)$ is a proper face of $\Gamma (S)$.

We define an equivalence relation in \ $\mathbb{R}_{+}^{n}$ by 
\begin{equation}
\begin{array}{ccc}
a\sim a^{\prime } & \text{iff} & F(a)=F(a^{\prime }).
\end{array}
\end{equation}
The equivalence classes of \ $\sim $ are sets of the form 
\begin{equation}
\Delta _{\tau }=\{a\in \mathbb{R}_{+}^{n}\mid F(a)=\tau \},
\end{equation}
where $\tau $ \ is a face of $\Gamma (S)$. We note that $\Delta _{\Gamma
(S)}=\{0\}$.

We recall that the cone strictly\ spanned \ by the vectors $%
a_{1},..,a_{e}\in \mathbb{R}_{+}^{n}\setminus \left\{ 0\right\} $ is the set 
$\Delta =$ \ $\left\{ \lambda _{1}a_{1}+...+\lambda _{e}a_{e}\mid \text{ }%
\lambda _{i}\in \mathbb{R}_{+}\text{, }\lambda _{i}>0\right\} $. If $%
a_{1},..,a_{e}$ are linearly independent over $\mathbb{R}$, $\Delta $ \ is
called \ a \textit{simplicial cone}. \ If moreover \ $a_{1},..,a_{e}\in 
\mathbb{Z}^{n}$, we say $\Delta $\ is \ a \textit{rational simplicial cone}.
If $\left\{ a_{1},..,a_{e}\right\} $ is a subset of a basis \ of the $%
\mathbb{Z}$-modulo $\mathbb{Z}^{n}$, we call $\Delta $ a \textit{simple cone}%
.

A precise description of the geometry of the equivalence classes modulo $%
\sim $ is as follows. Each facet \ $\gamma $ of $\Gamma (S)$\ has a unique \
primitive vector $a(\gamma )\in \mathbb{R}^{n}$ (i.e. a vector whose
components are positive integers and their greatest common divisor is one)
which is perpendicular to $\gamma $. We denote by $\frak{D}(\Gamma (S))$ the
set of such vectors. The equivalence classes are exactly the cones 
\begin{equation}
\Delta _{\tau }=\{\sum\limits_{i=1}^{e}\lambda _{i}a(\gamma _{i})\mid
\lambda _{i}\in \mathbb{R}\text{, }\lambda _{i}>0\},
\end{equation}
where $\tau $ runs through the set of faces of $\Gamma (S)$, and $\gamma
_{i} $, $i=1,..,e$\ are the facets containing $\tau $. We note that $\Delta
_{\tau }=\{0\}$ if and only if $\tau =\Gamma (S)$.

Each cone $\Delta _{\tau }$\ can be partitioned \ into a finite number of
simple cones $\Delta _{\tau ,i}$ . In addition, the subdivision can be
chosen such that each $\Delta _{\tau ,i}$ is spanned by part of $\frak{D}%
(\Gamma (S))$ (see e.g. \cite[sect. 2]{DH}). Thus from the above
considerations we have the following partition of $\mathbb{R}_{+}^{n}$:

\begin{equation}
\mathbb{R}_{+}^{n}=\{0\}\bigcup \bigcup\limits_{\tau \text{ }}\left(
\bigcup\limits_{i=1}^{l_{\tau }}\Delta _{\tau ,i}\ \right) ,  \label{5}
\end{equation}
where $\tau $ runs \ over the proper faces of $\Gamma (S)$, and each $\Delta
_{\tau ,i}$ \ is a simple cone contained in $\Delta _{\tau }$. \ We shall
say that $\left\{ \Delta _{\tau ,i}\right\} $ is a \textit{simple polyhedral
subdivision of }\ $\mathbb{R}_{+}^{n}$\ \textit{subordinated} to $\Gamma (S)$%
. \ 

\section{$p-$adic monomial integrals and Newton polyhedra}

We set $N_{i}=(N_{i,1},...,N_{i,n})\in \mathbb{N}^{n},N_{i}\neq 0,$ for $%
i=1,2,...,l,$ $\ N=(N_{1},...,N_{l})$, and $v\mathbf{=(}v_{1},...,v_{n}%
\mathbf{)}\in \mathbb{N}^{n}$. Let $x^{M_{i}}$ denote the monomial $%
\prod\limits_{j=1}^{n}x_{j}^{M_{i,j}}$, where $M_{i}=(M_{i,1},...,M_{i,n})%
\in \mathbb{N}^{n}$.

We associate to the pair $(N,v)$ the $p-$adic integral 
\begin{equation}
I_{(N,v)}(s):=\int\limits_{\left( \mathcal{P}_{K}^{e_{0}}\right)
^{n}}\left\| \left( c_{1}x^{N_{1}},...,c_{l}x^{N_{l}}\right) \right\|
_{K}^{s}\left| x^{v}\right| _{K}\left| dx\right| ,\text{ \ }Re(s)>0\text{,}
\label{6}
\end{equation}
where $c_{i}\in \mathcal{O}_{K}$, $i=1,2,...,l$, are constants, $e_{0}\in 
\mathbb{N}$, and \ $\left| dx\right| $\ is a Haar \ measure on $K^{n}$ so
normalized that $\mathcal{O}_{K}^{n}$\ has measure one. In this section we
show that \ $I_{(N,v)}(s)$ \ is a rational function of $q^{-s}$ and describe
its poles in terms \ of \ a Newton polyhedron.

We associate to $(N,v)$ the Newton polyhedron $\Gamma \left( N\right) .$ The
set of primitive perpendicular vectors to the faces of $\Gamma (N)$ is
denoted by $\frak{D}(\Gamma \left( N\right) )$. Each facet $\gamma \in
\Gamma \left( N\right) $ is the intersection of $\Gamma \left( N\right) $
and \ a \textit{supporting hyperplane} with equation \ $\left\langle
a(\gamma ),x\right\rangle =m(a(\gamma ))$. If $a(\gamma )=(a_{1},...,a_{n})$%
, we set $\sigma (a(\gamma )):=\sum\limits_{i=1}^{n}a_{i}$.

\begin{lemma}
\label{lemma}The integral \ $I_{(N,v)}(s)$ is\ a rational function of $%
q^{-s} $ with poles of the form 
\begin{equation}
s=-\frac{\sigma (a(\gamma ))+\left\langle v,a(\gamma )\right\rangle }{%
m(a(\gamma ))}+\frac{2\pi \sqrt{-1}k}{m(a(\gamma ))\log q},\text{ }k\in 
\mathbb{Z},  \label{7}
\end{equation}
where $\gamma $\ is a facet of $\Gamma (N),$ and $\left\langle a(\gamma
),x\right\rangle =m(a(\gamma ))$, with $m(a(\gamma ))$ $\neq 0$, \ is the
equation of the supporting hyperplane of $\gamma $.
\end{lemma}

\begin{proof}
First we fix some notation. Given $x=(x_{1},...,x_{n})\in \mathcal{O}%
_{K}^{n} $, we set $v(x):=(v(x_{1}),...,v(x_{n}))$, and 
\begin{equation}
E_{A}:=\{x\in \left( \mathcal{P}_{K}^{e_{0}}\right) ^{n}\mid v(x)\in A\cap 
\mathbb{N}^{n}\},  \label{8}
\end{equation}
for any $A\subseteq \mathbb{R}_{+}^{n}$.

We fix a simple polyhedral subdivision of \ $\mathbb{R}_{+}^{n}$\
subordinated to \ $\Gamma \left( N\right) $. From this subdivision we get
the following partition for \ $\left( \mathcal{P}_{K}^{e_{0}}\right) ^{n}$: 
\begin{equation}
\left( \mathcal{P}_{K}^{e_{0}}\right) ^{n}=\bigcup\limits_{\tau \subset
\Gamma \left( N\right) }\bigcup\limits_{i=1}^{l_{\tau }}E_{\Delta _{\tau ,i}}%
\text{, \ if \ \ }e_{0}\geqslant 1,  \label{9}
\end{equation}
where \ $\tau $ is a proper face of $\Gamma \left( N\right) $, and $\Delta
_{\tau ,i}$ \ is a simple cone. In the case in which \ $e_{0}$ $=0$, i.e. $%
\left( \mathcal{P}_{K}^{e_{0}}\right) ^{n}=\mathcal{O}_{K}^{n},$ \ the
following holds 
\begin{equation}
\mathcal{O}_{K}^{n}=\left( \mathcal{O}_{K}^{\times }\right) ^{n}\bigcup
\bigcup\limits_{\tau \subset \Gamma \left( N\right)
}\bigcup\limits_{i=1}^{l_{\tau }}E_{\Delta _{\tau ,i}}\text{.}  \label{10}
\end{equation}

From the above partitions, it follows that \ $I_{(N,v)}(s)$ \ is a finite
sum of integrals of the form 
\begin{equation}
\int\limits_{\left( \mathcal{O}_{K}^{\times }\right) ^{n}}\left\| \left(
c_{1}x^{N_{1}},...,c_{l}x^{N_{l}}\right) \right\| _{K}^{s}\left| dx\right|
=\left\| \left( c_{1},...,c_{l}\right) \right\| _{K}^{s}(1-q^{-1})^{n},\text{
\ }  \label{11}
\end{equation}
and \ 
\begin{equation}
\int\limits_{E_{\Delta _{\tau ,i}}}\left\| \left(
c_{1}x^{N_{1}},...,c_{l}x^{N_{l}}\right) \right\| _{K}^{s}\left|
x^{v}\right| _{K}\left| dx\right| \text{.}  \label{12}
\end{equation}
Without loss of generality, we assume that \ \ $e_{0}\geqslant 1$. By virtue
of the above considerations it is sufficient to show that integral (\ref{12}%
) is a rational function of $q^{-s}$ and that its poles have form (\ref{7}).

On the other hand, since $\Delta _{\tau ,i}$ is \ a simple cone, it is
spanned by some $a_{1},...,a_{e}\in \frak{D}(\Gamma \left( N\right) )$, $%
e=e_{\tau ,i}$, and $\{a_{1},...,a_{e}\}$ is a subset of a basis of the $%
\mathbb{Z}-$module $\mathbb{Z}^{n}$. Thus 
\begin{equation}
\Delta _{\tau ,i}\cap \mathbb{N}^{n}=\{\sum\limits_{i=1}^{e}a_{i}y_{i}\mid
y_{i}\in \mathbb{N\setminus \{}0\}\}.  \label{13}
\end{equation}
We set \ $I=\{1,2,...,e\}$, $J\subseteq I$, and for a \ fixed $\alpha
_{0}\in $ $\mathbb{N\setminus \{}0\}$,

\begin{equation*}
A_{J}:=\{\sum\limits_{i=1}^{e}a_{i}y_{i}\mid y_{i}\in \mathbb{N\setminus }%
\left\{ 0\right\} \text{, and }y_{i}\geqslant \alpha _{0}+1\Leftrightarrow
i\in J\},
\end{equation*}

\begin{equation*}
E_{A_{J}}=\{x\in \left( \mathcal{P}_{K}^{e_{0}}\right) ^{n}\mid v(x)\in
A_{J}\}.
\end{equation*}

We \ subdivide $\Delta _{\tau ,i}\cap \mathbb{N}^{n}$ as 
\begin{equation}
\Delta _{\tau ,i}\cap \mathbb{N}^{n}=\bigcup\limits_{J\subseteq I}A_{_{J}}.
\label{14}
\end{equation}
From partition (\ref{14}), it follows that \ 
\begin{equation}
\int\limits_{E_{\Delta _{\tau ,i}}}=\sum_{J\subseteq
I}\int\limits_{E_{A_{J}}}:=\sum_{J\subseteq I}I(s,E_{A_{J}}).
\end{equation}
The integral $I(s,E_{A_{\emptyset }})$ \ has an analytic \ continuation to
the complex plane as a polynomial in $q^{-s}$. Hence it is sufficient to
show that the integral $I(s,E_{A_{J}})$ is a rational function of $q^{-s}$,
and that its poles have form (\ref{7}), for $J$ $\neq \emptyset $.

The integral $I(s,E_{A_{J}})$ admits the following \ series expansion:

\begin{eqnarray}
I(s,E_{A_{J}})=\int\limits_{E_{A_{J}}}\left\| \left(
c_{1}x^{N_{1}},...,c_{l}x^{N_{l}}\right) \right\| _{K}^{s}\left|
x^{v}\right| _{K}\left| dx\right|  \notag \\
=\sum\limits_{k\in A_{J}}\int\limits_{\{x\in \left( \mathcal{P}%
_{K}^{e_{0}}\right) ^{n}\mid v(x)=k\}}\left\| \left(
c_{1}x^{N_{1}},...,c_{l}x^{N_{l}}\right) \right\| _{K}^{s}\left|
x^{v}\right| _{K}\left| dx\right| .  \label{19}
\end{eqnarray}
We write $k=(k_{1},...,k_{n})\in \mathbb{N}^{n}$, $\sigma \left( k\right)
:=k_{1}+...+k_{n}$, and 
\begin{equation}
x_{j}=\pi ^{k_{j}}u_{j},\text{ \ }u_{j}\in \mathcal{O}_{K}^{\times },\text{ }%
j=1,2,...,n\text{.}  \label{20}
\end{equation}
Then 
\begin{eqnarray}
\left| dx\right| =q^{-\sigma (k)}\left| du\right| ,  \label{21} \\
x^{N_{j}}=\pi ^{\left\langle k,N_{j}\right\rangle }u^{N_{j}}.  \notag
\end{eqnarray}
We know that $F(k)=\tau $. Thus 
\begin{equation}
\left\langle k,N_{j}\right\rangle =m(k),\text{ \ for every }N_{j}\in
\{N_{1},...,N_{l}\}\cap \tau \text{,}  \label{22}
\end{equation}
and 
\begin{equation}
\left\langle k,N_{j}\right\rangle >m(k)\text{, \ for every }N_{j}\in
\{N_{1},...,N_{l}\}\setminus \tau .  \label{23}
\end{equation}
We may assume, possibly after a renaming of the $N_{j},$ that $%
\{N_{1},...,N_{l}\}\cap \tau =\{N_{1},...,N_{r}\},$ with $r\leq l.$ Thus $%
\{N_{1},...,N_{l}\}\setminus \tau =\{N_{r+1},...,N_{l}\}$. With this
notation, and since $k=\sum\limits_{i=1}^{e}a_{i}y_{i}$, it follows from (%
\ref{21}), (\ref{22}), (\ref{23}), that 
\begin{equation}
x^{N_{j}}=\left\{ 
\begin{array}{c}
\pi ^{\sum\limits_{i=1}^{e}y_{i}m(a_{i})}u^{N_{j}},\text{ \ }j=1,2,...,r%
\text{,} \\ 
\pi ^{\sum\limits_{i=1}^{e}y_{i}\left\langle a_{i},N_{j}\right\rangle
}u^{N_{j}},\text{ \ \ }j=r+1,r+2,...,l\text{.}
\end{array}
\right. \text{ }  \label{24}
\end{equation}
In addition, for $j=r+1,2,...,l$, it holds that $\left\langle
a_{i},N_{j}\right\rangle >m(a_{i})$, for some index \ $i$.

We \ fix $\alpha _{0}$ \ satisfying the condition: 
\begin{equation*}
\alpha _{0}=e_{0}+\max_{1\leq i\leq n}\{v(c_{i})\}.
\end{equation*}

The constant $\alpha _{0}$ \ was selected to guarantee \ that the following
two conditions are satisfied: 
\begin{equation}
E_{A_{J}}=\left\{ x\in \left( \mathcal{P}_{K}^{e_{0}}\right) ^{n}\mid
v(x)=\sum\limits_{i=1}^{e}a_{i}y_{i}\text{, for some }y_{i}\in \mathbb{%
N\setminus }\left\{ 0\right\} \text{, }y_{i}\geqslant \alpha
_{0}+1\Leftrightarrow i\in J\text{ }\right\}  \label{24a}
\end{equation}
for every $J\neq \emptyset $, and 
\begin{equation}
\left\| \left( c_{1}x^{N_{1}},...,c_{l}x^{N_{l}}\right) \right\| _{K}\mid
_{E_{A_{J}}}=\left\| \left( c_{1}x^{N_{1}},...,c_{r}x^{N_{r}}\right)
\right\| _{K}\text{ }\mid _{E_{A_{J}}}.  \label{24b}
\end{equation}

Then from (\ref{21}), (\ref{24}), and (\ref{24b}) it follows that 
\begin{eqnarray}
\left\| \left( c_{1}x^{N_{1}},...,c_{l}x^{N_{l}}\right) \right\|
_{K}^{s}\left| x^{v}\right| _{K}\left| dx\right|  \label{25} \\
=b_{J}^{s}\text{ }q^{-\sum\limits_{i=1}^{e}y_{i}\left( \sigma
(a_{i})+\left\langle v,a_{i}\right\rangle +m(a_{i})s\right) }\left|
du\right| ,\text{ }  \notag
\end{eqnarray}
where $b_{J}$ is \ a constant, and $y_{i}\geqslant \alpha _{0}+1,$ $i\in J$.
Finally, (\ref{19}), (\ref{25}), and (\ref{24a}) imply that 
\begin{eqnarray}
I(s,E_{A_{J}})=\int\limits_{E_{A_{J}}}\left\| \left(
c_{1}x^{N_{1}},...,c_{l}x^{N_{l}}\right) \right\| _{K}^{s}\left|
x^{v}\right| _{K}\left| dx\right|  \notag \\
=b_{J}^{s}(1-q^{-1})^{n}\sum\limits_{
\begin{array}{c}
y_{i}\leq \alpha _{0}\text{ } \\ 
i\notin J
\end{array}
}\sum\limits_{
\begin{array}{c}
y_{i}\geqslant \alpha _{0}+1 \\ 
i\in J
\end{array}
}q^{-\sum\limits_{i=1}^{e}y_{i}\left( \sigma (a_{i})+\left\langle
v,a_{i}\right\rangle +m(a_{i})s\right) }  \notag \\
\text{ }=b_{J}^{s}(1-q^{-1})^{n}(Q_{J}(q^{-s}))\text{ }\prod\limits_{i\in
J}\left( \frac{q^{-\left( \alpha _{0}+1\right) \left( \sigma
(a_{i})+\left\langle v,a_{i}\right\rangle +m(a_{i})s\right) }}{1-q^{-\left(
\sigma (a_{i})+\left\langle v,a_{i}\right\rangle +m(a_{i})s\right) }}\right)
,  \label{26}
\end{eqnarray}

where $Q_{J}(q^{-s})$ is a polynomial in $q^{-s}$.
\end{proof}

\subsection{Remarks\label{remark}}

(1) For $l=1$ the \ previous lemma yields to a well-know result about $p-$%
adic elementary integrals (see e.g. \ \cite[Lemma 8.2.1]{I2}). (2) Given a
Newton polyhedron $\Gamma \left( N\right) \subseteq \mathbb{R}^{n}$, the
data $\sigma \left( a(\gamma )\right) $, $m(a(\gamma ))$, with $a(\gamma
)\in \frak{D}\left( \Gamma \left( N\right) \right) $, \ are invariant under
any renaming of the coordinates of $\mathbb{R}^{n}$.\ Thus by the
considerations made at the end of subsection 2.1, we can take any local
coordinate system in the computation of each $(N,v)$. (3) Lemma (\ref{lemma}%
) \ can be extended to integrals over sets of the form $b+\left( \mathcal{P}%
_{K}^{e_{0}}\right) ^{n}$, with $b\in \mathcal{O}_{K}^{n}$. Indeed, if\ $%
b\notin \left( \mathcal{P}_{K}^{e_{0}}\right) ^{n}$ \ and $\left|
x^{N_{i}}\right| _{K}\mid _{b+\left( \mathcal{P}_{K}^{e_{0}}\right) ^{n}}$
is constant for some $i$, the integral has an analytic continuation to the
complex plane as a polynomial in $q^{-s}$. In the other case, by a change of
variables, the integral can be reduced to an integral of type (\ref{6}).

\section{ Local zeta function for algebraic sets}

The first result of this paper is \ the following.

\begin{theorem}
\label{mainth}Let $f_{i}(x)\in K[x]$, $x=(x_{1},...,x_{n}),$ $%
f_{i}(0)=0,i=1,.,l$, be non-constant polynomials, and \ $f=(f_{1},...,f_{l})$%
. Let $\Phi $ be a Bruhat-Schwartz function. The local zeta function $%
Z_{\Phi }(s,f)$ admits a meromorphic continuation to the complex plane as a
rational function of $q^{-s}$ with poles in the set 
\begin{equation}
\bigcup\limits_{(N,v)}\bigcup\limits_{a\in \frak{D}(\Gamma \left( N\right)
)}\left\{ -\frac{\sigma (a)+\text{ }\left\langle v,a\right\rangle }{m(a)}+%
\frac{2\pi \sqrt{-1}k}{m(a)\log q},\text{ }k\in \mathbb{Z}\right\} ,
\label{30}
\end{equation}
where $(N,v)$ runs over the numerical data of $f$, and \ $\left\langle
a,x\right\rangle =m(a)$,$\ m(a)\neq 0$, is the supporting hyperplane \ of
the facet of \ $\Gamma \left( N\right) $ that corresponds to $a\in \frak{D}%
(\Gamma \left( N\right) )$.
\end{theorem}

\begin{proof}
By applying Hironaka's resolution theorem to the divisor $D_{K}=\cup
_{i=1}^{l}f_{i}^{-1}(0)$, \ we get \ $h:X_{K}\rightarrow A_{K}^{n}$ with $h$
a proper $K-$analytic map and $X_{K}$ a $n-$dimensional $K-$analytic
manifold. At every point $b$ of $X_{K}$ $\ $\ we can choose a chart $(U,\phi
_{U})$ such that $U$ contains $b$, $\phi _{U}(y)=(y_{1},...,y_{n})$ and \ 
\begin{eqnarray}
f_{i}\circ h &=&\epsilon _{i}\prod\limits_{j\in J}y_{j}^{N_{i,j}}\text{, }%
i=1,...,l,  \label{31} \\
h^{\ast }(\bigwedge\limits_{1\leq j\leq n}dx_{j}) &=&\left( \eta
\prod\limits_{j\in J}y_{j}^{v_{j}-1}\right) \bigwedge\limits_{1\leq j\leq
n}dy_{j},  \notag
\end{eqnarray}
with $J=\{j\in T\mid b\in E_{j}\}$, \ and $\epsilon _{i},\eta $ units of the
local ring $\mathcal{O}_{b}$ of \ $X_{K}$ at $b$. Since \ $h$ is proper and $%
A=Supp(\Phi )$ is compact open, we see that $B=h^{-1}(A)$ is compact open of 
$X_{K}$. Therefore, we can express $B$ \ as a finite disjoint union of
compact open sets$\ B_{\alpha }$ in $X_{K}$, \ satisfying $B_{\alpha
}\subseteq U_{\alpha }$. Since $\Phi $, $\left| \epsilon _{i}\right|
_{K},\left| \eta \right| _{K}$, are locally constant, \ after subdividing \ $%
B_{\alpha },$ we may assume that \ $\left( \Phi \circ h\right) \mid
_{B_{\alpha }}=\Phi \left( h(b)\right) $, $\left| \epsilon _{i}\right|
_{K}\mid _{B_{\alpha }}=\left| \epsilon _{i}(b)\right| _{K}$, $\left| \eta
\right| _{K}\mid _{B_{\alpha }}=\left| \eta (b)\right| _{K}$, and further $%
\phi _{U_{\alpha }}(B_{\alpha })=D_{\alpha }$, with $D_{\alpha }=w+\pi
^{e_{0}}\mathcal{O}_{K}^{n}$ for some $w\in K^{n}.$ Since $h:$ $X_{K}$ $%
\setminus h^{-1}\left( D_{K}\right) \rightarrow A_{K}^{n}$ $\setminus D_{K}$
is $K-$bianalytic, we then have 
\begin{equation}
Z_{\Phi }(s,f)=\sum\limits_{\alpha }\Phi \left( h(b)\right) \left| \eta
(b)\right| _{K}\int\limits_{^{D_{\alpha }}}\left\| \left( \left| \epsilon
_{1}(b)\right| _{K}y^{N_{1}},...,\left| \epsilon _{l}(b)\right|
_{K}y^{N_{l}}\right) \right\| _{K}^{s}\left| y^{v}\right| _{K}\left|
dy\right| ,  \label{32}
\end{equation}
where $N_{i}:=(N_{i,1},...,N_{i,n})\in \mathbb{N}^{n},i=1,...,l,$ $%
N:=(N_{1},...,N_{l}),v:=(v_{1}-1,...,v_{n}-1)\in \mathbb{N}^{n}$, and $%
N_{i,j}=0$, $v_{j}=1$, for $j$ not in $J$. \ The result follows \ by
applying lemma \ref{lemma}, and remark (\ref{remark}) (3) to the integral in
the right side of \ (\ref{32}).
\end{proof}

\bigskip

The proof of the main theorem is a generalization of Igusa's proof for the
case $l=1$, (see \cite[Theorem 8.2.1]{I2}).

\subsection{Remark}

Since the support of $\Phi $ is compact and $h$ is proper the set of
numerical data is finite.

\section{ Twisted local zeta function for algebraic sets}

Let $\chi _{i}$ be a character of $\mathcal{O}_{K}^{\times }$, i.e. \ a
homomorphism $\chi _{i}:\mathcal{O}_{K}^{\times }$\ $\rightarrow \mathbb{C}%
^{\times }$ with finite image, $i=1,...,l$. We formally put $\chi _{i}(0)$ $%
=0$, \ $i=1,...,l$, and $\chi :=(\chi _{1},...,\chi _{l})$. If $%
f=(f_{1},...,f_{l})$, we define 
\begin{equation*}
\chi (ac\text{ }f(x)):=\prod\limits_{i=1}^{l}\chi _{i}(ac\text{ }f_{i}(x)),
\end{equation*}
where \ $ac$ $z=z\pi ^{-v(z)}$ denotes the angular component of $z\in K$.
With the above notation, we associate to\ $\chi $ and $f$ the following
twisted local \ zeta function

\begin{equation}
Z_{\Phi }(s,\chi ,f):=\int\limits_{K^{n}}\Phi (x)\chi (ac\text{ }%
f(x))\left\| f(x)\right\| _{K}^{s}\mid dx\mid ,\,\,\,\,\,\,\,\,s\in \mathbb{C%
},\mathbb{\,}
\end{equation}
\noindent for $Re(s)>0$, where $\Phi $ is a Schwartz-Bruhat function, and $%
\mid dx\mid $ denotes the Haar measure on $K^{n}$so normalized that $%
\mathcal{O}_{K}^{n}$ has measure $1$. The proof of the following result \ is
a simple generalization of the proof of theorem \ref{mainth}.

\begin{theorem}
\label{theorem2}Let $f_{i}(x)\in K[x]$, $x=(x_{1},...,x_{n}),$ $%
f_{i}(0)=0,i=1,.,l$, be non-constant polynomials, \ $f=(f_{1},...,f_{l})$,
and $\chi =(\chi _{1},...,\chi _{l})$. Let $\Phi $ be a Bruhat-Schwartz
function. The local zeta function $Z_{\Phi }(s,\chi ,f)$ admits a
meromorphic continuation to the complex plane as a rational function of $%
q^{-s}$ with poles in the set 
\begin{equation}
\bigcup\limits_{(N,v)}\bigcup\limits_{a\in \frak{D}(\Gamma \left( N\right)
)}\left\{ -\frac{\sigma (a)+\text{ }\left\langle v,a\right\rangle }{m(a)}+%
\frac{2\pi \sqrt{-1}k}{m(a)\log q},\text{ }k\in \mathbb{Z}\right\} ,
\end{equation}
where $(N,v)=\left( \left\{ N_{1},...,N_{l}\right\} ,v\right) $, runs over
the numerical data of $f$, and \ $\left\langle a,x\right\rangle =m(a)$,$\
m(a)\neq 0$, is the supporting hyperplane \ of the facet of \ $\Gamma
_{(N,v)}$ that corresponds to $a\in \frak{D}(\Gamma \left( N\right) )$.
\end{theorem}

\section{Local zeta Functions For Monomial Algebraic Sets}

We fix $\left\{ N_{i}\in \mathbb{N}^{n}\mid N_{i}\neq 0,i=1,...,l\right\} $,
and denote by $\Gamma \left( \left\{ N_{1},...,N_{l}\right\} \right) $ its
Newton polyhedron. If $g(x)=\sum\limits_{l}a_{l}x^{l}\in K[x]$, $%
x=(x_{1},..,x_{m})$, is a non-constant polynomial satisfying $g(0)=0$, we
set supp$(g):=\{l\in \mathbb{N}^{m}\mid a_{l}\neq 0\}$.

\begin{definition}
Let $f_{i}(x)\in K[x]$, $x=(x_{1},...,x_{n}),$ $f_{i}(0)=0,i=1,...,l$, be
non-constant polynomials.\ The $K-$algebraic \ \ set 
\begin{equation*}
V_{K}(K)=V_{K}=\left\{ z\in K^{n}\mid f_{i}(z)=0,i=1,...,l\right\}
\end{equation*}
is called a monomial algebraic set \ if the following conditions hold:

(1) the polynomials $f_{i}(x)$\ have the form 
\begin{equation}
f_{i}(x)=c_{i}x^{N_{i}}+g_{i}(x),\text{ }i=1,...,l,  \label{33}
\end{equation}
with $c_{i}\in \mathcal{O}_{K}^{\times }$, \ $N_{i}\in \mathbb{N}^{n}$, $%
N_{i}\neq 0$, $g_{i}(x)\in \mathcal{O}_{K}[x],i=1,...,l$;

(2) any $m\in \cup _{i=1}^{l}$supp$(g_{i})$ belongs to the interior of \ $%
\Gamma \left( \left\{ N_{1},...,N_{l}\right\} \right) $, in the usual
topology of $\mathbb{R}^{n}$;

(3) the $K-$singular locus of $V_{K}$ is contained in \ 
\begin{equation*}
\bigcup_{i=1}^{l}\left\{ z\in K^{n}\mid z^{N_{i}}=0\right\} .
\end{equation*}
\end{definition}

We define \ $\Gamma _{V_{K}}:=\Gamma \left( \left\{ N_{1},...,N_{l}\right\}
\right) $ \ to be the Newton polyhedron of $V_{K}$. Also, we associate to $%
V_{K}$ the local zeta function 
\begin{equation}
Z(s,V_{K})=\int_{\mathcal{O}_{K}^{n}}\left\| f_{1}(x),...,f_{l}(x)\right\|
_{K}^{s}\mid dx\mid ,\,\,\,\,\,\,\,\func{Re}(\,s)>0.\mathbb{\,}  \label{34}
\end{equation}
The following theorem describes explicitly the meromorphic continuation of $%
Z(s,V_{K})$ in terms of $\Gamma _{V_{K}}$, its proof is a simple
modification of \ the proof of \ \ lemma \ \ref{lemma}.

\begin{theorem}
\label{th}Let $V_{K}$ a monomial $K-$algebraic set with Newton polyhedron $%
\Gamma _{V_{K}}$. Fix \ a simple polyhedral subdivision of $\mathbb{R}%
_{+}^{n}$: 
\begin{equation}
\mathbb{R}_{+}^{n}=\{0\}\bigcup \bigcup\limits_{\tau \text{ }}\left(
\bigcup\limits_{i=1}^{l_{\tau }}\Delta _{\tau ,i}\ \right) ,  \label{35}
\end{equation}
where $\tau $ is a proper face of $\ \Gamma _{V_{K}}$, and each $\Delta
_{\tau ,i}$ \ is a simple cone \ having \ $\tau $ as their first meet locus.
Then 
\begin{equation}
\text{\ \ }Z(s,V_{K})=(1-q^{-1})^{n}+\sum\limits_{\tau
}\sum\limits_{i=1}^{l_{\tau }}\left( \prod\limits_{j=1}^{e_{\tau ,i}}\frac{%
q^{-\left( \sigma (a_{j})+m(a_{j})s\right) }}{1-q^{-\left( \sigma
(a_{j})+m(a_{j})s\right) }}\right) ,  \label{36}
\end{equation}
where $\tau $ is a proper face of $\ \Gamma _{V_{K}}$, and $%
a_{1},...,a_{e_{\tau ,i}}$, are the generators of the cone $\Delta _{\tau
,i} $.
\end{theorem}

\end{document}